\documentclass{amsart}%
\usepackage{amssymb}
\usepackage{amsmath}
\usepackage{amsfonts}%
\usepackage{graphicx}
\providecommand{\U}[1]{\protect\rule{.1in}{.1in}}

\theoremstyle{plain}

\numberwithin{equation}{section}
\begin{document}
\title[The Pitchfork Bifurcation]{The Pitchfork Bifurcation }
\author{Indika Rajapakse}
\address{\noindent University of Michigan}
\email{indikar@umich.edu}
\author{Steve Smale}
\address{City University of Hong Kong, University of California, Berkeley}
\email{smale@cityu.edu.hk}
\thanks{We extend thanks to James Gimlett and Srikanta Kumar at Defense Advanced
Research Projects Agency for support and encouragement.\bigskip}
\date{8-28-2016}

\begin{abstract}
We give development of a new theory of the Pitchfork bifurcation, which
removes the perspective of the third derivative and a requirement of symmetry.

\end{abstract}
\maketitle

\noindent\textbf{Contents\medskip}

\noindent1. Introduction

\noindent2. Normal form

\noindent3. A relationship of our normal form to the main biological example

\noindent4. General pitchfork theorem

\noindent5. A biological perspective \bigskip

\noindent1. \textbf{Introduction\medskip}

\noindent The normal form for the pitchfork bifurcation is described usually
for one variable (see Guckenheimer and Holmes) by
\begin{equation}
\text{ \ }\frac{dx}{dt}=F(x)=\mu x-x^{3},\text{ \ }x\in%
\mathbb{R}
^{1},\ \mu\in%
\mathbb{R}
^{1}.\tag{1}%
\end{equation}

\noindent Note that this equation is invariant under the change of the
variable $x\rightarrow-x.$ That is,$\ F$ is an odd function. This condition
suggests that the Pitchfork bifurcation is generic for problems that have
symmetry. To obtain the above form one argues (or assumes) that the second
derivative of $F(x)$ is zero.\ This normal form is standard in the literature
on Pitchfork bifurcation [2, 3]. \medskip

\noindent Our own proof [4] for the one variable case is different from the
previous literature in that a new uniformity condition is satisfied in place
of symmetry or the hypothesis, vanishing of a second derivative. For the new
uniformity condition see below. In addition we take the path of using this
uniformity together with the Poincar\'{e} Hopf Theorem to show that the second
derivative must be zero.\medskip

\noindent We have found little in the literature on the case of more than one
variable, except for suggesting that (1) still applies. Kuznetsov [5] has a
proof for an $n$ variable case for pitchfork that assumes invariance under an
involution. Kuznetsov's hypothesis eliminates a second derivative.\medskip

\noindent Our paper [4] states an $n-$dimensional version of the pitchfork
theorem. That proof involves a reduction to one dimensional theory using
center manifold theory. In this paper we give a proof, which gives new insight
especially for more than one variable.\medskip

\noindent As we will see the one dimensional case with emphasis on the third
derivative is misleading. In fact for two variables there is a robust
pitchfork bifurcation in a quadratic system (with no third derivative at all).
Our previous work [4] is a background for this paper.\medskip

\noindent\textbf{2. Normal form\medskip}

\noindent We propose that the \textbf{normal form} for the Pitchfork
bifurcation for the dimension of the space greater than one be:%

\begin{align}
\text{ \ }\frac{dx}{dt}  & =y^{2}-ay-x\tag{2}\\
\frac{dy}{dt}  & =x^{2}-ax-y,\text{ \ \ \ }x,y\geq0\nonumber
\end{align}

\noindent For each point of the bifurcation parameter $a,$ the "central
equilibrium" is $\left(  x,y\right)  =\left(  0,0\right)  ,$ and may be
described in terms of the two isoclines$-$the curves in the $(x,y)$ plane
where $\frac{dx}{dt}=0$ (the $x$ isoclines) and where $\frac{dy}{dt}=0$ (the
$y$ isoclines). \ For given $a$ the equilibria are given as the intersection
of these isoclines. The isoclines are described respectively, by the
equations:%
\begin{align}
x  & =y^{2}-ay\tag{3a}\\
y  & =x^{2}-ax.\tag{3b}%
\end{align}
Figure 1 shows the curves for selected values of $a.$ For $a\leq1,$ the
central equilibrium is the only equilibrium relevant for the bifurcation. For
$a>1$ there are four points of intersection and three equilibrium points for
the model. The transition between these situations occur at $a=1,$ when the
isoclines are tangent to each other at $(0,0)$.%
\begin{figure}[ptb]%
\centering
\includegraphics[
natheight=2.021100in,
natwidth=6.437700in,
height=1.3413in,
width=4.2125in
]%
{../../Users/indikar/Desktop/Pichf.tif}%
\caption{Geometry of isoclines for the normal form (3). Left: a single
equilibrium at $(1,1)$ for $a<1$; middle: a single equilibrium at $(1,1)$ for
$a=1$; right: for $a>1$, $(1,1)$ becomes unstable, and two stable equilibria
emerge flanking $(1,1)$. To view how changing the parameter $a$ changes the
isoclines, see: https://www.desmos.com/calculator/x5fz6lavrx }%
\end{figure}
\medskip

\noindent\textbf{The defining characteristic of a pitchfork bifurcation is the
transition from a single stable equilibrium to two new stable equilibria
separated by a saddle. \ The saddle emerges from the old stable equilibrium.}
\ \qquad\qquad\qquad\qquad\qquad\qquad\qquad\qquad\qquad\medskip

\noindent The Figure 1 illustrates this characteristic as $a$ increases from
less than one to greater than one.\medskip

\noindent Next we will obtain the intersection points of two isoclines as in
Equations 3a and 3b. The curves described by Equations 3a and 3b are quadratic
and one can solve their intersection points analytically. If we substitute 3b
into 3a, solving for $x$ yields for the solutions $x_{i}$ as can be checked:
\begin{align}
x_{1}  & =0,\text{ central equilibrium}\tag{4a}\\
x_{2}  & =\frac{1}{2}\left(  a-1\right)  +\frac{1}{2}\sqrt{\left(  a-1\right)
\left(  a+3\right)  }\nonumber\\
x_{3}  & =\frac{1}{2}\left(  a-1\right)  -\frac{1}{2}\sqrt{\left(  a-1\right)
\left(  a+3\right)  }\nonumber\\
x_{4}  & =a+1\nonumber
\end{align}
Similarly, substituting Equation 3a into 3b and solving for $y$ yields the
solutions $y_{i}$:
\begin{align}
y_{1}  & =0,\text{ central equilibrium}\tag{4b}\\
y_{2}  & =\frac{1}{2}\left(  a-1\right)  -\frac{1}{2}\sqrt{\left(  a-1\right)
\left(  a+3\right)  }\nonumber\\
y_{3}  & =\frac{1}{2}\left(  a-1\right)  +\frac{1}{2}\sqrt{\left(  a-1\right)
\left(  a+3\right)  }\nonumber\\
y_{4}  & =a+1.\nonumber
\end{align}

\noindent The pairs $\left(  x_{i},y_{i}\right)  ,$ $i=1,..,4$ describe the
equilibria. Note that the intersection point $\left(  x_{4},y_{4}\right)  $ is
extraneous to the bifurcation phenomena. Note also that the Equation 4a and 4b
show the bifurcation effect at $a=1$ and $a>1.$\medskip

\noindent The two isoclines are parabolas and they get translated vertically
and horizontally as the parameter $a$ increases. One can see the intersections
of these parabolas in terms of simple analytic geometry, and these
intersections include the equilibria of the pitchfork\textbf{.} This formalism
allows us to see the pitchfork variables in terms of geometry and extend the
analysis to the nonsymmetric case, \ $\frac{dx}{dt}=y^{2}-ay-x,$\ $\frac
{dy}{dt}=x^{2}-bx-y$.\medskip\ 

\noindent The Jacobian matrix $J$ of the first partial derivatives of System
(2) at the central equilibrium $\left(  x_{1},y_{1}\right)  =\left(
0,0\right)  $ of Equation 3 is:%
\begin{equation}
J=\left(
\begin{array}
[c]{cc}%
-1 & -a\\
-a & -1
\end{array}
\right)  .\tag{5}%
\end{equation}
For each $a$, the eigenvalues are $\lambda_{1}=a-1,\lambda_{2}=-a-1$ and the
corresponding eigenvectors are $\left(  -1,1\right)  $ and $(1,1)$
respectively. When $a=1,$ $\lambda_{1}=0$ and $\lambda_{2}<0.$\ When $a<1,$
both $\lambda_{1},\lambda_{2}$ have negative real parts. Hence the central
equilibrium is stable. \ When $a>1,$\ $\lambda_{1}>0$ and $\lambda_{2}<0,$ and
the central equilibrium is a saddle. The qualitative structure is
robust.\medskip

\noindent Eigenvalues at the equilibrium $\left(  x_{2},y_{2}\right)  $ are
given by: $\lambda_{1}=-1+\sqrt{-(a-1)(a+3)+3}$ and $\lambda_{2}%
=-1-\sqrt{-(a-1)(a+3)+3}.$ When $1<a,$ both $\lambda_{1},\lambda_{2}$ have
negative real parts. Hence the equilibrium is stable, similarly for $\left(
x_{3},y_{3}\right)  $. \medskip

\noindent\textbf{Remark 1}: When $1.\,\allowbreak2361<a,$ both $\lambda
_{1},\lambda_{2}$ are complex conjugate numbers with negative real parts. The
pitchfork phenomena continues after $a=1.2361$ using the equilibrium Equations
4a and 4b. \ \bigskip

\noindent\textbf{3. A relationship of our normal form (2) to a main biological
example}\medskip

\noindent We are motivated by work by Gardner et al. [6] for a "synthetic,
bistable gene-regulatory network \ldots\ [to] provide a simple theory that
predicts the conditions necessary for bistability." \ The toggle, as designed
and constructed by Gardner et al., is a network of two mutually inhibitory
genes that acts as a switch by some mechanism, as a control for switching from
one basin to another. Consider the particular setting of Gardner et al.'s
circuit design of the toggle switch as:%
\begin{align}
\frac{dx}{dt} &  =\frac{\alpha_{1}}{1+y^{m}}-x\tag{6}\\
\frac{dy}{dt} &  =\frac{\alpha_{2}}{1+x^{n}}-y.\nonumber
\end{align}
If $m=n=0,$ the equilibrium is $x=\frac{\alpha_{1}}{2},$ $y=\frac{\alpha_{2}%
}{2},$ and the eigenvalues of the Jacobian are negative. If $\alpha_{1}%
<2\max(x)$ and $\alpha_{2}<2\max(y),$ the system has a unique global stable
equilibrium [6, 7].\medskip\ 

\noindent When $\alpha_{1}=\alpha_{2}=2,$ and $m=n>0,$ the system in Equation
6 can be written more specifically as in Gardner et al.:
\begin{align}
\frac{dx}{dt} &  =\frac{2}{1+y^{m}}-x=f\left(  x,y\right) \tag{7}\\
\frac{dy}{dt} &  =\frac{2}{1+x^{m}}-y=g(x,y),\text{ \ \ \ }0\leq x,y.\nonumber
\end{align}
\ \ For this system with $0\leq m\leq2,$ every equilibria must be $\left(
1,1\right)  .$ \noindent We give the Taylor approximation for each $m$ about
the equilibrium $(1,1)$. For this approximation consider the derivatives
$f_{x},$ $f_{y},f_{xx},f_{yy},f_{xy}$ and similar for $g$, all at the point
$(1,1).$These can be computed as:
\begin{align*}
f_{x}  & =-1,\text{ }f_{y}=-\frac{1}{2}m,\text{ }g_{x}=-\frac{1}{2}m,\text{
}g_{y}=-1\\
f_{yy}  & =g_{xx}=\frac{1}{2}m,\text{ and the remaining derivatives are zero.}%
\end{align*}
\ \noindent Therefore, the Taylor approximation about the equilibrium $(1,1)$
for each $m\leq2$ is given (deleting the remainder):%
\begin{align}
\frac{dx}{dt}  & =\frac{1}{2}my^{2}-\frac{3}{2}my+m+1-x\tag{8a}\\
\frac{dy}{dt}  & =\frac{1}{2}mx^{2}-\frac{3}{2}mx+m+1-y\tag{8b}%
\end{align}
Recall our normal form Equation 1 written in terms of $u$ and $v$ is:%
\begin{align}
\frac{du}{dt}  & =v^{2}-av-u\tag{9a}\\
\frac{dv}{dt}  & =u^{2}-au-v.\tag{9b}%
\end{align}
We wish to compare Equation 8 and 9 with the parameter values at the
bifurcation points, namely $m=2$ for Equation 8 and $a=1$ for Equation 9. Then
the transformation $x=u+1$ and $y=v+1$ gives a correspondence between
Equations 8 and 9 at these parameter values. After the bifurcation points $m$
and $a$ vary dependently but we don't know the functional relationship. \ The
point is that each $m$ and $a$ increasing from the bifurcation value create a
pitchfork bifurcation.\bigskip

\noindent\textbf{4. General pitchfork theorem}

Recall some setting from our previous paper [4].%
\begin{equation}
\frac{dx}{dt}=F_{\mu}\left(  x\right)  ,\ \mu\in%
\mathbb{R}
,\text{ }x\in%
\mathbb{R}
^{n}\text{and }\left\vert \mu\right\vert ,\left\vert x\right\vert
<\varepsilon\tag{10}%
\end{equation}

\noindent This is associated to a family $F_{\mu}$ with bifurcation parameter
$\mu\in(-\varepsilon,\varepsilon)$ describing $\frac{dx}{dt}=F_{\mu}\left(
x\right)  .$ Here $x$ belongs to a domain $X$ of\textbf{\ }$%
\mathbb{R}
^{n},\ $and $F_{0}\left(  x\right)  =F\left(  x\right)  .$ We suppose that the
dynamics of $F_{\mu}$ is that of a stable equilibrium $x_{\mu},$ basin
$B_{\mu}$ for $\mu<\mu_{0},$ and that the bifurcation is at $\mu
_{0}.\smallskip$ \ We assume that$\ x_{\mu}=0$ is an equilibrium for all
$\mu.\medskip$

\noindent\textbf{Uniformity condition:}\ "First bifurcation from a stable
equilibrium." \noindent The equilibrium does not "leave it's basin" in the
sense that there is a neighborhood $N$ of $x_{0}$, such that $N$ is contained
in $B_{\mu}$ for all $\mu<0.$ \medskip

\noindent Note this implies by a uniform continuity argument, that even at
$\mu=0$, $\ x_{0}$ is a sink in the sense that $x\left(  t\right)  \rightarrow
x_{0},$ if the initial point belongs to $N.$ The dynamics of $F_{\mu},$
$\mu=0$ has a "basin" $B:$ $x_{0}$ is a "weak sink." It follows that in this
space $B,$ the only equilibria of $F_{\mu}$ is $x_{\mu}$ for all $\mu\leq0.$
One could say that $\mu_{0}$ is the "first" bifurcation.\medskip

\noindent Define $J_{\mu}$ to be the matrix of partial derivatives of $F_{\mu
}$ \ at $x_{\mu}$. \ The eigenvalues of $J_{\mu}$ for $\mu<\mu_{0} $ all have
negative real part, either real, or in complex conjugate pairs. At the
bifurcation, one has either a single real eigenvalue becoming zero and then
positive with the pitchfork (if det($J_{\mu_{0}})\neq0$ then a complex
conjugate pair of distinct eigenvalues with real parts zero becomes positive
after the bifurcation and then the Hopf oscillation occurs, not discussed
here). \medskip

\noindent A pitchfork bifurcation converts a stable equilibrium into two
stable equilibria (the Hopf bifurcation converts a stable equilibrium into a
stable periodic solution).\medskip

\noindent\textbf{The pitchfork bifurcation theorem:} \ In Equation 10, let
$x_{\mu}$ be a stable equilibrium for $\mu$ for all $\mu<\mu_{0}.$ Suppose the
uniformity condition is satisfied. If the determinant of $J_{\mu_{0}}=0,$ then
generically the dynamics undergoes a pitchfork bifurcation.\medskip

\noindent1. If $n=1,$ this has been proved in our paper [4] as discussed above.

\noindent2. For $n>1,$ generically there is a pitchfork as exemplified by the
normal form in Section 2. \ \medskip

\noindent\textbf{Sketched of proof of the pitchfork bifurcation theorem for
}$n>1$.\medskip

\noindent Note in Equation 10, the equation for the equilibria is $F_{\mu
}\left(  x\right)  =0$ for each $\mu.\medskip$

\noindent Consider the hypothesis above and consider the equilibrium $x_{\mu}%
$. At $\mu=\mu_{0}$, the determinant of $J_{\mu_{0}}$ becomes zero. Then by
the local stable manifold theory, this equilibrium changes from a sink
$\mu<\mu_{0}$, to a saddle for $\mu>\mu_{0}.$ This saddle has an $n-1$
dimensional contracting stable manifold $\mathbf{W}_{\mu}^{u}$ and a one
dimensional expanding stable manifold $\mathbf{W}_{\mu}^{s}.$\medskip

\noindent We will be using the following.\medskip

\noindent\textbf{Poincar\'{e}-Hopf index theorem}: Suppose $\frac{dx}%
{dt}=F\left(  x\right)  $, $x$ belongs to $X$ and $F:X\rightarrow%
\mathbb{R}
^{n}.$ Suppose $X$ \ homeomorphic to a closed ball [8] and $F\left(  x\right)
$ points to the interior of $X$ for each $x$ belonging to the boundary
\begin{equation}
\sum_{F(x)=0,\text{ }x\in X}sign\left(  \det(J)\right)  \left(  x\right)
=\left(  -1\right)  ^{n}\tag{11}%
\end{equation}
where $J$\textbf{\ }is the Jacobian of $F$ at $x.$ The formula on the left
side of Equation (11) is the Poincar\'{e}-Hopf index$.$\ In particular
generically, in the case $n$ is odd, there is an odd number of equilibria.
\medskip

\noindent Then we observe that the Poincar\'{e}-Hopf index at the equilibrium
for the sink is $\left(  -1\right)  ^{n}$ for $\mu<\mu_{0}$. \ For $\mu
>\mu_{0}$ of this equilibria changes to $\left(  -1\right)  ^{n+1}$ as the
sink changes to a saddle. This follows from the eigenvalue structure of the
saddle. \ Therefore, by the Poincar\'{e}-Hopf index theorem there must be
other new equilibria for $\mu>\mu_{0}$. Generically these new equilibria must
be two in number and they are sinks. We have obtained the defining property of
the pitchfork.\medskip

\noindent The above needs to be carried out uniformly for all $\mu.$ This
procedure follows a suggestion of Mike Shub [9, 10]. \medskip

\noindent We consider the 2-dimensional center manifold, $C_{M}$, at $\mu
=\mu_{0}$ and the equilibrium $x=x_{0}$ with the added equation $\frac{d\mu
}{dt}=0.$ This gives a dynamical system of $n+1$ equations. The $C_{M}$ is two
dimensional and the stable manifold $\mathbf{W}_{\mu}^{u}$ is $n-1\dim
$ensional, corresponding to the eigenvalues $J_{\mu_{0}}$ with negative real
parts. \ The $C_{M}$ corresponds to the span of the null space (eigenvector
corresponding the zero eigenvalue) and the space of the variable $\mu.$ The
$C_{M}$ projects on to $\mu$ and the inverse image of $\mu$ is the expanding
one dimensional manifold of the equilibrium $x_{0},$ for $F_{\mu}$.\bigskip

\noindent\textbf{5. A biological perspective} \medskip

\noindent In our previous work [4], a cell can be thought as a point in the
basin and the cell type can be identified with a basin. Thus, the identity of
a specific cell type in our genome dynamics can be defined by characteristic
gene expression pattern at the equilibrium. \ We suggest that the emergence of
a new cell type from this original cell type, through differentiation,
reprogramming, or cancer a result of pitchfork bifurcation, is a departure
from a stable equilibrium and requires cell division. In normal cell division
during differentiation or reprogramming, a cell can undergo symmetric or
asymmetric division. Let A be the mother cell in the following. In symmetric
division, two identical daughter cells arise: A and A (B and B), that have
genomes with the same activity [11, 12]. In asymmetric division, two daughter
cells arise: A and B, that have genomes with different activity. Both cases
above reflect a pitchfork bifurcation. In another type of asymmetric division,
two daughter cells arise, B and C, where the activity of both genomes is
different from the mother also reflecting a pitchfork bifurcation. One example
of this is in cases of abnormal cell division, where chromosomes are
mis-segregated, resulting in one daughter with too many and one daughter with
too few chromosomes. This type of division may be one of the initiating events
in emerging cancer cells [13]. Capturing these events in terms of our
bifurcations may give us insight into emergence of a cell type.\bigskip

\noindent\textbf{Acknowledgments: }We would like to thank Lindsey Muir, Thomas
Ried, and especially Mike Shub for helpful discussions.


\bigskip

\begin{thebibliography}{99}                                                                                               %
\bibitem[1]{}Guckenheimer J, Holmes P (2002) Nonlinear Oscillations, Dynamical
Systems, and Bifurcations of Vector Fields (Springer, New York).

\bibitem[2]{}Strogatz S (2000) Non-linear Dynamics and Chaos: With
applications to Physics, Biology, Chemistry and Engineering (Perseus Books).

\bibitem[3]{}Wiggins S (2003) Introduction to applied nonlinear dynamical
systems and chaos (Springer-Verlag, New York, second edition).

\bibitem[4]{}Rajapakse I, Smale S (2016) Mathematics of the genome.
\textit{Foundations of Computational Mathematics}, 1-23.

\bibitem[5]{}Kuznetsov Y (1998) Elements of applied bifurcation theory
(Springer-Verlag, New York, second edition).

\bibitem[6]{}Gardner T, Cantor C, Collins J (2000) \ Construction of a genetic
toggle switch in Escherichia coli. \textit{Nature,} 403: 339-342.

\bibitem[7]{}Ellner S, Guckenheimer J (2006) Dynamics models in biology.
(Princeton University Press, Princeton, New Jersey).

\bibitem[8]{}V. Guillemin, A. Pollack. Differential Topology. Prentice-Hall, 1974

\bibitem[9]{}Shub M (1987) Global Stability of Dynamical Systems
(Springer-Verlag, New York).

\bibitem[10]{}Hirsch M, Pugh C, Shub M (1977) Invariant Manifolds, Lecture
Notes in Mathematics, Vol. 583 (Springer-Verlag, New York).

\bibitem[11]{}Hartwell H, Kastan M (1994) Cell cycle control and cancer.
\textit{Science}, 266 (5192): 1821-1828.

\bibitem[12]{}Massagu\'{e} J (2004) G1 cell-cycle control and cancer
\textit{Nature} 18; 432 (7015): 298-306.

\bibitem[13]{}Ried T et. al. (2012) The consequences of chromosomal aneuploidy
on the transcriptome of cancer cells. \textit{Biochim Biophys Acta.} 1819(7): 784-93.
\end{thebibliography}
\end{document}